\documentclass[letterpaper,11pt,oneside]{amsart}

\usepackage{amssymb}
\usepackage{amsxtra}
\usepackage[mathscr]{eucal}
\usepackage{graphics}

\numberwithin{equation}{section}

\newtheorem{theorem}{Theorem}[section]
\newtheorem{proposition}[theorem]{Proposition}
\newtheorem{lemma}[theorem]{Lemma}

\newtheorem{conjecture}[theorem]{Conjecture}

\theoremstyle{definition}

\newtheorem{example}[theorem]{Example}

\newtheorem{remark}[theorem]{Remark}

\newcommand{\al}{\alpha}

\newcommand{\artp}{\al_{R_2^{(p)}/\fp}}
\newcommand{\be}{\mathbf{e}}
\newcommand{\bn}{\mathbf{n}}
\newcommand{\fp}{\mathfrak{p}}
\newcommand{\Newt}{\mathscr{N}}

\newcommand{\ZZ}{\mathbb{Z}}
\newcommand{\RR}{\mathbb{R}}

\begin{document}

\allowdisplaybreaks

\title[Asymptotic geometry of
non-mixing sequences]{Asymptotic geometry of\\non-mixing sequences}
\author{Manfred Einsiedler and Thomas Ward}

\address{ME,
Mathematisches Institut,
 Universit\"at Wien, Strudlhofgasse 4, A-1090, Vienna, Austria;
Mathematics Department,
State College, PA. 16802, U.S.A.}
\email{manfred@mat.univie.ac.at}

\address{TW, School of Mathematics, University of East
  Anglia, Norwich NR4 7TJ, United Kingdom}
\email{t.ward@uea.ac.uk}

\thanks{The first author was supported
by EPSRC Grant GR/M 49588.}

\date{\today}


\begin{abstract}
The exact order of mixing for zero-dimensional algebraic dynamical
systems is not entirely understood. Here non-Archimedean norms in
function fields of positive characteristic are used to exhibit an
asymptotic shape in non-mixing sequences for algebraic
$\ZZ^2$-actions. This gives a relationship between the order of
mixing and the convex hull of the defining polynomial. Using these
methods, we show that an algebraic dynamical system for which any
shape of cardinality three is mixing is mixing of order three, and
for any $k\ge1$ exhibit examples that are $k$-fold mixing but not
$(k+1)$-fold mixing.
\end{abstract}

\maketitle

\section{Introduction and definitions}

Let $\alpha$ be a $\ZZ^d$-action by $\mu$-preserving
transformations of a non-trivial probability space
$(X,\mathcal S,\mu)$.
A sequence
$({\mathbf n}_1^{(j)},{\mathbf n}_2^{(j)},\dots,{\mathbf n}_r^{(j)})$
of $r$-tuples of elements of $\ZZ^d$ with
\begin{equation}\label{movingapart}
{\mathbf n}_s^{(j)}-{\mathbf n}_t^{(j)}\to\infty
\mbox{ as }
j\to\infty\mbox{ for every }s\neq t.
\end{equation}
is {\sl mixing for} $\alpha$ if for any sets
$A_1,\dots,A_r\in\mathcal S$,
\begin{equation}\label{mixingrsequence}
\lim_{j\to\infty}
\mu\left(
\alpha^{-{\mathbf n}_1^{(j)}}(A_1)\cap\dots\cap
\alpha^{-{\mathbf n}_r^{(j)}}(A_r)\right)=
\mu(A_1)\cdots\mu(A_r).
\end{equation}
If any sequence satisfying (\ref{movingapart}) has
(\ref{mixingrsequence}) then $\alpha$ is {\sl mixing of order
$r$}. As usual, an action that is mixing of order 2 is called
simply {\sl mixing}. The maximum value of $r$ for which
(\ref{movingapart}) implies (\ref{mixingrsequence}) -- if this is
finite -- is the {\sl order of mixing} of $\alpha$. If
(\ref{movingapart}) implies (\ref{mixingrsequence}) for all $r$,
then $\alpha$ is {\sl mixing of all orders}. A finite set
$\{{\mathbf n}_1,\dots,{\mathbf n}_r\}$ of integer vectors is a
{\sl mixing shape} for $\alpha$ if
\begin{equation}
\lim_{k\to\infty} \mu\left( \alpha^{-k{\mathbf
n}_1}(A_1)\cap\dots\cap \alpha^{-k{\mathbf
n}_r}(A_r)\right)=\mu(A_1)\cdots\mu(A_r).
\end{equation}
In general, mixing properties of {\sl shapes} have no bearing on
the {\sl order} of mixing. For example, there are non-mixing
actions for which all shapes are mixing (see \cite[Theorem
1.2]{MR99e:28031}). An algebraic dynamical system -- one in which
$X$ is assumed to be a compact metrizable abelian group, $\mu$ is
Haar measure, and $\alpha$ acts by automorphisms -- with all
shapes mixing is, in contrast, mixing of all orders. However, it
is not clear whether non-mixing shapes detect the exact order of
mixing for algebraic systems that are not mixing of all orders --
see \cite[Sect. 27, 28]{MR97c:28041} for an overview of this
problem, and \cite{MR80b:28030} for Ledrappier's seminal example
which showed that mixing algebraic $\ZZ^2$-actions need not be
mixing of all orders.

\begin{conjecture}\label{conjr}An algebraic dynamical system for
which all shapes of cardinality $r$ are mixing is mixing of order
$r$.
\end{conjecture}

As remarked above, the (suitably interpreted) conjecture holds for
``$r=\infty$'' in the sense that all shapes mixing implies mixing
of all orders for algebraic dynamical systems. It also holds when
$r=2$ -- that is, if each element $\alpha^{\mathbf n}$, ${\mathbf
n}\neq0$ is mixing, then the whole action $\alpha$ is mixing (see
\cite[Theorem 1.6]{MR97c:28041}). In Theorem \ref{conj3} below we
show that the conjecture holds when $r=3$ and $d=2$. Moreover, a
weaker lower bound for the order of mixing for $d=2$ is shown in
Theorem \ref{mixingfromcoding}, and this is used to give examples
with any given order of mixing in Section \ref{sec:example}.

\section{Algebraic $\ZZ^2$-actions}

Let $R_2=\ZZ[u_1^{\pm1},u_2^{\pm1}]$, the ring of Laurent
polynomials with integer coefficients in the commuting variables
$u_1$, $u_2$. Following \cite{MR91g:22008}, associate a given
algebraic $\ZZ^2$-action $\alpha$ on $X$ to an $R_2$-module as
follows. Let $M$ be the countable Pontryagin (character) group of
$X$, and define $u_j\cdot m=\widehat{\al}^{\be_j}(m)$ for all
$m\in M$, where $\be_j\in\ZZ^2$ is the $j$th unit vector and
$\widehat{\al}^{\be_j}$ is the automorphism of $M$ dual to
~$\al^{\be_j}$ for $j=1,2$. A polynomial $F\in R_2$ has the form
$F(u)= \sum_{\bn\in\ZZ^2}c_F(\bn)u^{\bn}$, where $c_F(\bn)\in\ZZ$
and $c_F(\bn)=0$ for all but finitely many $\bn\in\ZZ^2$, and
$u^{\bn}=u_1^{n_1}u_2^{n_2}$. Then $F\cdot m=\sum_{\bn\in\ZZ^2}
c_F(\bn)\widehat{\al}^{\bn}(m)$ for every $m\in M$. Conversely,
suppose that $M$ is a countable $R_2$-module, and let
$X_M=\widehat{M}$ be its compact abelian character group. Each
$u_j$ is a unit in $R_2$, so the map $\gamma_j$ defined by
$\gamma_j(m)=u_j\cdot m$ is an automorphism of $M$. Define an
algebraic $\ZZ^2$-action $\al_M$ on $X_M$ by
$\al_M^{\be_j}=\widehat{\gamma}_j$. Thus
$(X_M,\al_M)\leftrightarrow M$ gives a one-to-one correspondence
between algebraic $\ZZ^2$-actions and countable $R_2$ modules. See
\cite[Chap.~II]{MR97c:28041} for further explanation and many
examples.

Approximating indicator functions of the sets $A_1,\dots,A_r$
appearing in (\ref{mixingrsequence}) with trigonometric
polynomials shows that for the algebraic action $\alpha$
corresponding to the $R_2$-module $M$, (\ref{mixingrsequence}) is
equivalent to the statement that for any $m_1,\dots,m_r\in M$, not
all zero, the relation
\begin{equation}\label{algform}
m_1u^{{\mathbf n}_1^{(j)}}+ m_2u^{{\mathbf n}_2^{(j)}}+\dots+
m_ru^{{\mathbf n}_r^{(j)}}=0
\end{equation}
holds for only finitely many $j$.

In order to use this to make progress with Conjecture \ref{conjr},
the following well-known results are needed. A consequence of the
algebraic formulation of mixing in \eqref{algform} is that a given
sequence $({\mathbf n}_1^{(j)},{\mathbf n}_2^{(j)},\dots,{\mathbf
n}_r^{(j)})$ is mixing for $\al_M$ on $X_M$ if and only if it is
mixing for $\al_{R_2/\mathfrak P}$ on $X_{R_2/\mathfrak P}$ for
every prime ideal $\mathfrak P\subset R_2$ associated to the
module $M$ (see \cite[Theorem~3.3]{MR95f:28022} or
\cite[Theorem~2.2]{MR95c:22011} for a proof of this). Moreover, if
$X_M$ is connected (equivalently, if $M$ is torsion-free as an
additive group), then by \cite{MR95c:22011}, if $\al_M$ is mixing,
it is also mixing of all orders. If $\mathfrak P$ is a prime ideal
generated by a single prime number then $\al_{R_2/\mathfrak P}$ is
the $\ZZ^2$-shift on $\{0,\ldots,p-1\}^{\ZZ^2}$ and is therefore
mixing of all orders. Finally, if $\mathfrak P$ contains a prime
number and is not of the form $p\cdot R_2+F\cdot R_2$ for some
non-zero polynomial $F$, then $R_2/\mathfrak P$ is finite and
$\al_{R_2/\mathfrak P}$ is not mixing. It follows that in order to
show Conjecture \ref{conjr} for $d=2$, it is sufficient
to study algebraic dynamical systems corresponding to cyclic
modules of the form $R_2/\mathfrak P$ for some prime ideal
$\mathfrak P\subset R_2$ of the form $\langle p,F\rangle=p\cdot
R_2+F\cdot R_2$, where $p$ is a prime number and $F\in R_2$.

\section{Asymptotic geometry}

From now on we will only be working in rings of characteristic $p$
for a fixed prime $p$, so we replace the coefficient ring $\ZZ$
with the finite field $\mathbb F_p=\ZZ/p\ZZ$. Notationally this
replaces the polynomial $F\in R_2$ with the reduced mod $p$
polynomial
$f=\overline{F}\in R_2^{(p)}=\mathbb F_p[u_1^{\pm1},u_2^{\pm2}]$,
and the prime ideal
$$\mathfrak P=<p,F>=p\cdot R_2+F\cdot R_2$$
with the prime
ideal
$$\mathfrak p=<f>=f\cdot R_2^{(p)}.$$
Notice that $R_2/\mathfrak P$ is then the same ring as
$R_2^{(p)}/\mathfrak p$. For simplicity of notation we will write
$\alpha=\artp$ for the $\ZZ^2$-action on the dual of
${R_2^{(p)}/\mathfrak p}$ where $\mathfrak p=\langle f\rangle$ for
a fixed irreducible polynomial $f\in R_2^{(p)}$.

Let $\mathscr M(\alpha)$ denote the {\sl order of mixing} -- the
largest value of $r$ for which (\ref{movingapart}) implies
(\ref{mixingrsequence}) -- of $\alpha$. Finding $\mathscr
M(\alpha)$ is difficult (see \cite[Sect.~28]{MR97c:28041}) even
for this special class of systems. Our purpose here is to find a
new inequality that relates $\mathscr M(\alpha)$ to the {\sl
shape} of the polynomial $f$.

The polynomial $f$ can be written $f(u)=
\sum_{\bn\in\ZZ^2}c_f(\bn)u^{\bn}$, where $c_f(\bn)\in\mathbb F_p$
for all $\bn\in\ZZ^2$ and $c_f(\bn)=0$ for all but finitely many
$\bn\in\ZZ^2$. Let
$$S(f)=\{\bn\in\ZZ^2\mid c_f(\bn)\neq0\}$$
denote the {\sl support} of $f$, and
$\Newt(f)$ the
{\sl convex hull} of $S(f)$. Thus $S(f)$ is some
finite set of points in $\ZZ^2$, and $\Newt(f)$ is a convex
polygon in $\ZZ^2$.

The relationships found below between faces of $\Newt(f)$ and
measurable properties of the action
are essentially equivalent to
the geometry of half-spaces and relative entropies in the paper
\cite{MR95f:28022}. The approach taken here using non-Archimedean
norms seems to be better adapted to attacking Conjecture
\ref{conjr}.

\begin{theorem}\label{mixingfromcoding}
Assume that $\Newt(f)$ is an $R$-gon and $f$ is irreducible. Let
$\alpha$ be the $\ZZ^2$-action on the dual of $R_2^{(p)}/\langle
f\rangle$. Then
$$R-1\le\mathscr M(\alpha)
<\vert S(f)\vert.$$
\end{theorem}

\begin{theorem}\label{conj3}
Conjecture {\rm\ref{conjr}} holds when $r=3$ for
algebraic $\mathbb{Z}^2$-actions.
\end{theorem}

The method of proof of Theorem~\ref{mixingfromcoding} is to show
that an arbitrary non-mixing sequence for $\alpha$ must
asymptotically reflect part of the structure of $\Newt(f)$ (the
slopes of the faces). Theorem~\ref{conj3} then holds because the
slopes of a triangle determine its shape.

The key step in the proof is to construct norms that reflect the
geometry of $\Newt(f)$. To clarify this, an example is described
in detail (see Example \ref{exampleofvaluation}). Before doing
this, some background on non-Archimedean norms and Newton polygons
is needed (see \cite[Chap.~2]{goss} or
\cite[Chap.~IV.3]{koblitz} for more details).
A {\em non-Archimedean norm} $\vert\cdot\vert$ on an integral
domain $S$ is a function
\[
 \vert\cdot\vert:S\rightarrow \RR
\]
with the properties
\begin{enumerate}
\item[(i)] $\vert a\vert\geq 0$,
\item[(ii)] $\vert a\vert=0$ only for $a=0$,
\item[(iii)] $\vert a+b\vert\leq \max(\vert a\vert,\vert b\vert)$, and
\item[(iv)] $\vert ab\vert=\vert a\vert\cdot\vert b\vert$
\end{enumerate}
for every $a,b\in S$. A norm always extends uniquely to the field $K$
of quotients of $S$ (or to any intermediate ring between
$S$ and $K$). An immediate
consequence of property (iii) is that given any equation
in $S$ of the form
\begin{equation}\label{lateron}
 \sum_{k=1}^na_k=0\mbox{ with }a_k\neq 0,
\end{equation}
there must be at least two indices $i\neq j$ for which
\begin{equation}\label{alongtheshore}
 \max_{1\le k\le n}\vert a_k\vert =\vert a_i\vert =\vert a_j\vert.
\end{equation}
The {\sl Newton polygon} of a polynomial
$g(x)=g_0+g_1x+\dots+g_nx^n\in S[x]$
with respect to the non-Archimedean norm $\vert\cdot\vert$
is defined to be the highest convex polygonal line joining
$(0,-\log\vert g_0\vert)$ to $(n,-\log\vert g_n\vert)$ which
passes on or below all the points $(j,-\log\vert g_j\vert)$ for
$j=0,\dots,n$. The {\sl vertices} of the Newton polygon are the points
$(j,-\log\vert g_j\vert)$ where the slope changes, and the slope
of a line in the Newton polygon joining the vertices
$(i,-\log\vert g_i\vert)$ and $(j,-\log\vert g_j\vert)$ is
$(\log\vert g_j\vert-\log\vert g_i\vert) /(i-j)$. The basic
property of the Newton polygon is the following: for each slope
$\lambda$, there is a non-Archimedean norm
$\vert\cdot\vert^{(\lambda)}$ on the extension ring $S[x]/\langle g\rangle$,
coinciding with $\vert\cdot\vert$ on the constant polynomials
(identified with $S$), and with $\vert
x\vert^{(\lambda)}=\exp(\lambda)$.

Up to a natural equivalence the (non-trivial) non-Archimedean
norms on $S=\mathbb F_p[u^{\pm1}]$ are those of the form
\begin{equation}\label{e:gprime}
 \vert f\vert _g=p^{-\operatorname{ord}_gf},
\end{equation}
(where $g\in\mathbb F_p[u]$ is an irreducible polynomial, and
$\operatorname{ord}_gf$ denotes the multiplicity of $g$ in the
prime decomposition of $f$), together with the exceptional norm
\[
 \vert f\vert _\infty=p^{\operatorname{deg}f}.
\]

\begin{example}\label{exampleofvaluation}
Let $f(u_1,u_2)=u_2+u_1+u_1^3u_2$, and view $f$ as an element of
${\mathbb F}_p[u_2^{\pm1}][u_1^{\pm1}]$ (cf. Figure \ref{fig1},
showing the support of $f$ as dots).
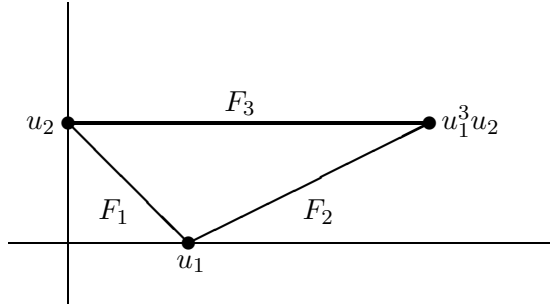
\begin{figure}[h]
\setlength{\unitlength}{0.08cm}
\begin{picture}(80,60)
\put(-10,10){\line(1,0){90}} \put(0,0){\line(0,1){50}} \thicklines
\put(0,30){\line(1,-1){20}} \put(20,10){\line(2,1){40}}
\put(0,30){\line(1,0){60}} \put(5,14){$F_1$} \put(39,14){$F_2$}
\put(26,32){$F_3$} \put(18,6){$u_1$} \put(-7,29){$u_2$}
\put(62,29){$u_1^3u_2$}
\put(20,10){\circle*{2}}
\put(0,30){\circle*{2}}
\put(60,30){\circle*{2}}
\end{picture}
\caption{\label{fig1}The faces of $\Newt(u_2+u_1+u_1^3u_2)$}
\end{figure}
Choose the norm $\vert\cdot\vert=\vert \cdot\vert _{u_2}$ on
${\mathbb F}_p[u_2^{\pm1}]$, so that $\vert
u_2\vert=\frac{1}{p}$.

How this norm extends to the ring extension $R={\mathbb
F}_p[u_2^{\pm1}][u_1^{\pm1}]/\langle f\rangle$ is determined by
the Newton polygon of $f$ viewed as a polynomial for $u_1$ with
coefficients in ${\mathbb F}_p[u_2^{\pm1}]$:
\begin{equation}\label{ringpoly}
f(u_1)=u_2\cdot u_1^0 + 1\cdot u_1^1+0\cdot u_1^2+u_2\cdot u_1^3.
\end{equation}
The four points that define the Newton polygon (Figure \ref{fig2})
are therefore
$$(0,-\log_p\vert u_2\vert)=(0,1),$$
$$(1,-\log_p\vert 1\vert)=(1,0),$$
$$(2,-\log_p\vert0\vert)=(2,\infty)$$and
$$(3,-\log_p\vert u_2\vert)=(3,1)$$(logarithms base $p$ are used
for convenience, and $\vert 0\vert=0$). Notice that the Newton
polygon of $f$ shown in Figure \ref{fig2} does not coincide with
the convex hull of the support in Figure \ref{fig1} -- but they
have the same faces pointing towards negative powers of $u_2$ in
Figure \ref{fig1} (equivalently, towards monomials for which
$\vert\cdot\vert$ is big).

\begin{figure}[h]
\setlength{\unitlength}{0.08cm}
\begin{picture}(70,40)
\put(0,10){\line(1,0){70}} \put(0,10){\line(0,1){30}} \thicklines
\put(0,30){\line(1,-1){20}} \put(20,10){\line(2,1){40}}
\put(-4,28){$1$} \put(19,4){$1$} \put(39,4){$2$} \put(59,4){$3$}
\put(0,30){\circle*{2}} \put(20,10){\circle*{2}}
\put(60,30){\circle*{2}}
\end{picture}
\caption{\label{fig2}The Newton polygon of $f\in{\mathbb
F}_p[u_2^{\pm1}][u_1^{\pm1}]$ with respect to $\vert\cdot\vert$}
\end{figure}
From the Newton polygon in Figure \ref{fig2} it follows that there
are two norms $\vert\cdot\vert_1$, $\vert\cdot\vert_2$ extending
$\vert\cdot\vert$ to $R$; the first has $\vert
u_1\vert_1=\frac{1}{p}$ (from the line segment with slope $-1$)
and the second has $\vert u_1\vert_2=\root\of{p}$ (from the line
segment with slope $1/2$).

Thus the vector $\left(\begin{matrix} \log_p\vert
u_1\vert_1\\\log_p\vert u_2\vert_1\end{matrix}\right)=
\left(\begin{matrix}-1\\-1\end{matrix}\right)$
is an outward normal to the face $F_1$ of $\Newt(f)$.
The same expression
using $\vert\cdot\vert_2$ gives $\left(\begin{matrix} \log_p\vert
u_1\vert_2\\\log_p\vert u_2\vert_2\end{matrix}\right)=
\left(\begin{matrix}1/2\\-1\end{matrix}\right)$, an outward normal
to the face $F_2$.

If the norm $\vert\cdot\vert^{\prime}=\vert\cdot\vert_\infty$ on
${\mathbb F}_p[u_2]$ is chosen initially, then $\vert
u_2\vert^{\prime}=p$ (so that the monomials with big norm are in
the upper-half plane of Figure \ref{fig1}) and the corresponding
Newton polygon is determined from \eqref{ringpoly} by the points
$$
(0,-\log_p\vert u_2\vert)=(0,-1),$$
$$(1,-\log_p\vert 1\vert)=(1,0),$$
$$(2,-\log_p\vert0\vert)=(2,\infty)$$and
$$(3,-\log_p\vert u_2\vert)=(3,-1).$$
The resulting Newton polygon is shown in Figure \ref{fig2andabit};
it shows there is only one extension to $R$, and the resulting
norm $\vert\cdot\vert_1^{\prime}$ has the property that
$\left(\begin{matrix} \log_p\vert u_1\vert_1^{\prime}\\\log_p\vert
u_2\vert_1^{\prime}\end{matrix}\right)=
\left(\begin{matrix}0\\1\end{matrix}\right)$ is an outward normal
to the face
$F_3$ of $\Newt(f)$.
\begin{figure}[h]
\setlength{\unitlength}{0.08cm}
\begin{picture}(70,40)(0,-25)
\put(0,10){\line(1,0){70}} \put(0,10){\line(0,-1){30}}\thicklines
\put(0,-10){\circle*{2}}\put(20,10){\circle*{2}}\put(60,-10){\circle*{2}}
\put(0,-10){\line(1,0){60}}\put(-8,-11){$-1$} \put(19,12){$1$}
\put(39,12){$2$} \put(59,12){$3$}
\end{picture}
\caption{\label{fig2andabit}The Newton polygon of $f\in{\mathbb
F}_p[u_2^{\pm1}][u_1^{\pm1}]$ with respect to
$\vert\cdot\vert^{\prime}$}
\end{figure}
\end{example}

Example \ref{exampleofvaluation} generalizes to the next lemma.

\begin{lemma}\label{lemma}
For each face $F$ of the convex hull $\Newt(f)$, there is a norm
$\vert\cdot\vert^{(F)}$ on the ring $R_2^{(p)}/\langle f\rangle$
with the property that the vector $\left(\begin{matrix}
\log_p\vert u_1\vert^{(F)}\\\log_p\vert
u_2\vert^{(F)}\end{matrix}\right)$ is an outward normal to
the face $F$ of $\Newt(f)$.
\end{lemma}

\begin{proof}
Choose a face $F$ of $\Newt(f)$. If necessary, exchange $u_1$ and
$u_2$ so that $F$ is not a vertical line. By replacing $u_2$ with
$u_2^{-1}$ if necessary, assume further that $F$ is one of the
{\sl lower} faces of $\Newt(f)$ when drawn in the plane (that is,
like $F_1$ or $F_2$ in Figure \ref{fig1}). Let
$\vert\cdot\vert=\vert\cdot\vert_{u_2}$ be the norm corresponding
to the irreducible polynomial $g=u_2$ in $\mathbb F_p[u_2]$ as in
(\ref{e:gprime}). We may also assume (after multiplying $f$ by a
suitable monomial) that $f$ is a polynomial in $u_1$ with
coefficients in $\mathbb F_p[u_2]$,
 \begin{equation}\label{e:fqi}
  f=\sum_{k=0}^nq_i(u_2)u_1^k\mbox{ with }q_0q_n\neq 0.
 \end{equation}
For each non-zero coefficient $q_i$, let $m_i$ be the largest $m$
for which $u_2^{m}$ divides $q_i(u_2)$ in $\mathbb F_p[u_2]$. By the
definition of the norm $\vert\cdot\vert$,
\[
  -\log_p\vert q_i(u_2)\vert=m_i.
 \]
On the other hand, the coefficient of $u_2^{m_i}$ in $q_i$ is
nonzero, so $(i,m_i)\in S(f)$. This shows that the points
$(i,m_i)$ appear both in the {\sl Newton polygon} of $f$
(considered as a polynomial in $u_1$ as in (\ref{e:fqi})) and in
the {\sl support} of $f$. Comparing $\Newt(f)$ and the Newton
polygon shows that the lower faces of $\Newt(f)$ comprise exactly
the Newton polygon of $f$. Thus, there is a norm
$\vert\cdot\vert^{(F)}$ on $R_2^{(p)}/\langle
 f\rangle$ extending $\vert\cdot\vert$ for which
 \[
  \vert u_1\vert^{(F)}=p^\lambda,\mbox{ where }\lambda\mbox{ is the slope of $F$.}
 \]
The vector
 \[
\left(\begin{matrix} \log_p\vert u_1\vert^{(F)}\\\log_p\vert
u_2\vert^{(F)}\end{matrix}\right)=\left(\begin{matrix} \lambda\\
-1
\end{matrix}\right)
 \]
is therefore an outward normal to the face $F$.
\end{proof}

These norms will now be used to show that a non-mixing sequence
must asymptotically approximate some of the structure of
$\Newt(f)$, by applying the simple
observation that (\ref{lateron})
implies (\ref{alongtheshore}) to the algebraic characterisation of
mixing, using a norm adapted to the shape $\Newt(f)$.

\begin{proposition}\label{nastyprop}
Assume that $\left(A^{(j)}\right)= ({\mathbf n}_1^{(j)},{\mathbf
n}_2^{(j)},\dots,{\mathbf n}_r^{(j)})$ is a sequence in
$\left(\ZZ^2\right)^r$ with the property that
\begin{equation}\label{nonmixeq}
m_1u^{{\mathbf n}_1^{(j)}}+ m_2u^{{\mathbf n}_2^{(j)}}+\dots+
m_ru^{{\mathbf n}_r^{(j)}}=0
\end{equation}
for all $j$, where $m_1,\dots,m_r\in R_2^{(p)}/\fp$ are non-zero.
Write $\Newt(A^{(j)})$ for the convex hull of the set $\{{\mathbf
n}_1^{(j)},{\mathbf n}_2^{(j)},\dots,{\mathbf n}_r^{(j)}\}$. Fix a
face $F$ of $\Newt(f)$. Then there is a constant $K>0$ such that
there is a face of $\Newt(A^{(j)})$ spanned (without loss of
generality) by ${\mathbf n}_1^{(j)},{\mathbf n}_2^{(j)}$, and a
vector ${\mathbf m}^{(j)}_{\vphantom{i}}$ with the property that
the line through ${\mathbf n}_1^{(j)}, {\mathbf
m}^{(j)}_{\vphantom{i}}$ is parallel to $F$ and $\Vert{\mathbf
m}^{(j)}-{\mathbf n}_2^{(j)}\Vert\le K$.
\end{proposition}

\begin{proof}
Pick a face $F$ of $\Newt(f)$ and use Lemma \ref{lemma} to find a
norm $\vert\cdot\vert$ on $R_2^{(p)}/\langle f\rangle$ so that
$\left(\begin{matrix} \log_p\vert u_1\vert\\\log_p\vert
u_2\vert\end{matrix}\right)$ is an outward normal to $\Newt(f)$
through $F$. Let $L=2\max_{i=1,\dots,r}\{ \left\vert{\log_p\vert
m_i\vert}\right\vert\}.$

Fix $j$ for now, and choose $t$ to maximize $\vert u^{\mathbf
n^{(j)}_t}\vert$. This corresponds to $\mathbf n^{(j)}_t$ being
extremal in $\Newt(A^{(j)})$ in the direction of
$\left(\begin{matrix} \log_p\vert u_1\vert\\\log_p\vert
u_2\vert\end{matrix}\right)$. Let $\ell$ be the line through
$\mathbf n^{(j)}_t$  parallel to $F$. Assume that no other point
of $A^{(j)}$ lies within $L$ of the line $\ell$. Then for $i\neq
t$ we get
\[
\mathbf n^{(j)}_i \left(\begin{matrix} \log_p\vert
u_1\vert\\\log_p\vert u_2\vert\end{matrix}\right) < \mathbf
n^{(j)}_t \left(\begin{matrix} \log_p\vert u_1\vert\\\log_p\vert
u_2\vert\end{matrix}\right)-L
\]
and
\[
 \vert m_iu^{\mathbf
n^{(j)}_i}\vert=\vert m_i\vert\vert u^{\mathbf n^{(j)}_i}\vert<
p^{-L}\vert m_i\vert \vert u^{\mathbf n^{(j)}_t}\vert\leq \vert
m_tu^{\mathbf n^{(j)}_t}\vert.
\]
This shows that in (\ref{nonmixeq}) one term is bigger than all
the others, which is impossible for a non-Archimedean norm
$\vert\cdot\vert$. Therefore, for every $j$ there must be a second
point $n^{(j)}_s$ within distance $L$ of the line $\ell$. If
necessary, pass to a subsequence so that $t$ and $s$ are
independent of $j$. By renaming the indices assume $t=1,s=2$.

Let $V$ be the rational subspace normal to $\left(\begin{matrix}
\log_p\vert u_1\vert\\\log_p\vert u_2\vert\end{matrix}\right)$,
and choose a basis $\mathbf b_1,\mathbf b_2$ of $\ZZ^2$ with
$\mathbf b_1\in V$. For every $\mathbf n\in\ZZ^2$ the component
$a_2\mathbf b_2$ of $\mathbf n=a_1\mathbf b_1+a_2\mathbf b_2$ can
be found by projection along $V$. Since $\mathbf n_2^{(j)}$ is
within $L$ of the line $\ell=V+n_1^{(j)}$, the projections of
$\mathbf n_1^{(j)}$ and $\mathbf n_2^{(j)}$ onto $\mathbf b_2\ZZ$
are close together, say within distance $K$. This shows that
$\mathbf n_2^{(j)}=\mathbf m^{(j)}+\mathbf c^{(j)}$ with $\mathbf
m^{(j)}\in V+n_1^{(j)}$, $\parallel \mathbf c^{(j)}\parallel < K$.
\end{proof}

\begin{proof} (of Theorem \ref{mixingfromcoding})
First recall that $S(f)$ is automatically a non-mixing shape for
$\alpha$ (see \cite[Examples 27.1]{MR97c:28041}), so
$$\mathscr M(\alpha)
<\vert S(f)\vert.$$ On the other hand, the convex hull
$\Newt(A^{(j)})$ is a convex polygon. By Proposition
\ref{nastyprop}, each face $F$ of $\Newt(f)$ must appear with a
uniformly bounded error as one of the faces of $\Newt(A^{(j)})$.
Since the differences in (\ref{movingapart}) go to infinity, the
slope of this matching face approaches the slope of $F$. It
follows that there must be at least $R$ faces on the convex set
$\Newt(A^{(j)})$, so $A^{(j)}$ must have at least $R$ points. Thus
$R-1\leq \mathscr M(\alpha)$.
\end{proof}

\begin{proof} (of Theorem \ref{conj3})
If $\Newt(f)$ lies on a line, then $\alpha$ cannot be mixing. If
$\Newt(f)$ is an $R$-gon with $R>3$ then Theorem
\ref{mixingfromcoding} shows that $\mathscr M(\alpha)\ge3$. So
assume that $\Newt(f)$ is a triangle, with vertices $\mathbf d_i$
for $i=1,2,3$. Let $F_i$ be the face of $\Newt(f)$ spanned by
$\mathbf d_i,\mathbf d_{i+1}$ (reduce subscripts mod $3$). Assume
additionally that $\alpha$ is {\sl not} mixing of order $3$: We
will deduce that there is a {\sl shape} of cardinality $3$ which
is not mixing.

There exist non-zero elements $a,b,c\in R_2^{(p)}/\langle
f\rangle$, and three sequences $\mathbf n_i^{(j)}$ of lattice
points with \eqref{movingapart} and
\begin{equation}\label{e:abc}
au^{{\mathbf n}_1^{(j)}}+ bu^{{\mathbf n}_2^{(j)}}+ cu^{{\mathbf
n}_r^{(j)}}=0
\end{equation}
for all $j$. Without loss of generality, assume that $\mathbf
n_i^{(j)},\mathbf n_{i+1}^{(j)}$ (again, reduce subscripts mod
$3$) are the two points in $A^{(j)}$ satisfying Proposition
\ref{nastyprop} for the face $F_i$. Let $i=1$. Then there exists
$\mathbf m^{(j)}\in\ZZ^2$ bounded away from $\mathbf n_2^{(j)}$
such that $\mathbf m^{(j)}$ is actually on the line through
$\mathbf n_1^{(j)}$ parallel to $F_1$. Passing to a subsequence,
we may assume that
\[
  \mathbf n_2^{(j)}-\mathbf m^{(j)}=\mathbf k
\]
is independent of $j$. Transform equation (\ref{e:abc}) to
\[
 au^{{\mathbf n}_1^{(j)}}+ (bu^{\mathbf k})u^{{\mathbf m}^{(j)}}+ cu^{{\mathbf
n}_r^{(j)}}=0.
\]
In other words, by changing $b$ slightly we can assume that
${\mathbf n}_1^{(j)}$ and ${\mathbf n}_2^{(j)}$ are actually on a
line parallel to $F$.
 By changing $\mathbf n_2^{(j)}$
if necessary by another bounded amount, we can assume that there
exists $s_j\in\mathbb N$ with
\[
 \mathbf n_2^{(j)}-\mathbf n_1^{(j)}=s_j(\mathbf d_2-\mathbf d_1).
\]
Repeat this for the face $F_2$, and change the vector ${\mathbf
n}_3^{(j)}$ accordingly. Now we know that there exists
$t_j\in\mathbb N$ such that
\[
 \mathbf n_3^{(j)}-\mathbf n_2^{(j)}=t_j(\mathbf d_3-\mathbf d_2).
\]
Using Proposition  \ref{nastyprop} again for the face $F_3$, we
see that $ \mathbf n_1^{(j)}$ and $ \mathbf n_3^{(j)}$ are almost
on the same line parallel to $F_3$. This almost similarity between
$\Newt(f)$ and $\Newt(A^{(j)})$ shows that $|s_j-t_j|$ is bounded.
Changing $\bn^{(j)}_3$, again by a uniformly bounded amount, we
can achieve $s_j=t_j$ so the triangles are similar. This shows
that $\{\mathbf d_1,\mathbf d_2,\mathbf d_3\}$ is a non-mixing
shape.
\end{proof}

\section{Examples}\label{sec:example}

\begin{example}
Theorem \ref{mixingfromcoding} shows that if $f$ is an irreducible
polynomial for which the support $S(f)$ coincides with the extreme
points of the Newton polygon $\Newt(f)$, then $\mathscr
M(\alpha)=\vert S(f)\vert-1$. In order to produce an example
$\alpha$ with prescribed order of mixing $\mathscr M(\alpha)=k$,
it is therefore sufficient to exhibit such an irreducible
polynomial with $\vert S(f)\vert=k+1$. This may be done using
Eisenstein's irreducibility criterion (see \cite{MR87i:11172} for
a general norm-theoretic treatment of the Eisenstein criterion).
Two simple examples will illustrate the method; it is clear from
these how to build an example for any order of mixing. We are
grateful to Klaus Schmidt for pointing out that an explicit
construction of such a family of examples was not known
previously.
\begin{enumerate}
\item To find an example with order of mixing $3$,
consider $f(u_1,u_2)=u_1^2+u_1u_2^2+u_2^3+u_2\in \mathbb
F[u_2][u_1]$; the prime $u_2\in\mathbb F[u_2]$ divides the
coefficients $u_2^2$ and $u_2^3+u_2$ but $u_2^2$ does not divide
the coefficient $u_2^3+u_2$. By Eisenstein's criterion $f$ is
irreducible. The support of the polynomial is shown in Figure
\ref{fig3}.
\begin{figure}[h]
\setlength{\unitlength}{0.06cm}
\begin{picture}(80,70)
\put(0,0){\line(1,0){65}}
\put(0,0){\line(0,1){65}}
\thicklines
\put(0,20){\line(2,-1){40}}
\put(0,20){\line(0,1){40}}
\put(0,60){\line(1,-1){20}}
\put(20,40){\line(1,-2){20}}
\put(38,-6){$u_1^2$}
\put(-8,19){$u_2$}
\put(-8,59){$u_2^3$}
\put(22,39){$u_1u_2^2$}
\put(0,20){\circle*{2}}
\put(0,60){\circle*{2}}
\put(40,0){\circle*{2}}
\put(20,40){\circle*{2}}
\end{picture}
\caption{\label{fig3}The support of a polynomial giving $3$-fold mixing}
\end{figure}
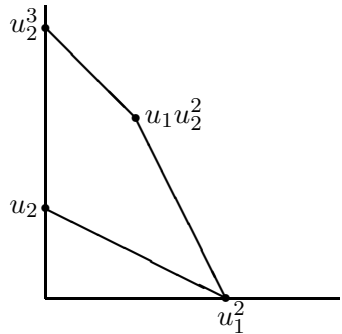
\item To find an example with order of mixing $5$,
let $f(u_1,u_2)=u_1^6+u_1^5u_2+u_1^3u_2^2+u_2+u_2^3.$
As before, this is seen to be irreducible by viewing
it as a polynomial in $u_1$ with coefficients in
$\mathbb F[u_2]$.
The support of the polynomial is shown in Figure \ref{fig4}.
\begin{figure}[h]
\setlength{\unitlength}{0.06cm}
\begin{picture}(135,70)
\put(0,0){\line(1,0){130}}
\put(0,0){\line(0,1){65}}
\thicklines
\put(0,20){\line(6,-1){120}}
\put(0,20){\line(0,1){40}}
\put(0,60){\line(3,-1){60}}
\put(60,40){\line(2,-1){40}}
\put(100,20){\line(1,-1){20}}
\put(118,-6){$u_1^6$}
\put(120,0){\circle*{2}}
\put(-8,19){$u_2$}
\put(0,20){\circle*{2}}
\put(-8,59){$u_2^3$}
\put(0,60){\circle*{2}}
\put(60,40){\circle*{2}}
\put(62,40){$u_1^3u_2^2$}
\put(100,20){\circle*{2}}
\put(102,19){$u_1^5u_2$}
\end{picture}
\caption{\label{fig4}The support of a polynomial giving $4$-fold mixing}
\end{figure}
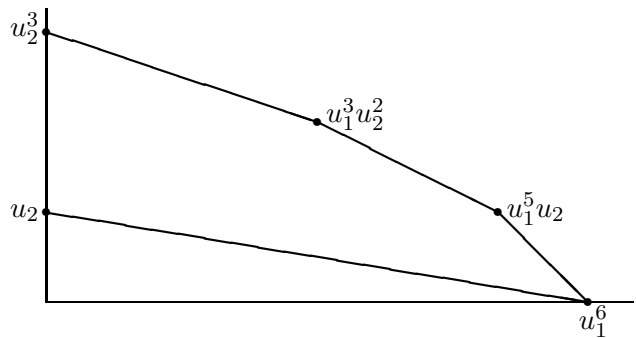
\end{enumerate}
Notice that in these examples we are
choosing the shape of the support freely;
it is also possible to find examples for which
any prescribed shape is the minimal non-mixing shape by
\cite{MR99e:28031}, though not in a constructive fashion.
\end{example}

\begin{example}
Theorem \ref{conj3} shows that the system corresponding
to the ideal $\fp=\langle 2, 1+u_1+u_2+u_2^2\rangle$
is 3-mixing, answering a question in \cite[p.~283]{MR97c:28041}.
\end{example}

\begin{example}In the previous example, we used
the fact from \cite{MR97c:28041} that no shape
with cardinality 3 is non-mixing.
An alternative method to show this is to use
a result of Voloch on solutions to
$ax+by=1$ in function fields.
Consider again $\fp=\langle 2, 1+u_1+u_2+u_2^2\rangle$;
then Theorem \ref{mixingfromcoding} says that
$$2\le\mathscr M(\alpha)
<4,$$ and we wish to show that $\mathscr M(\alpha)=3$. To see
this, assume that
$$\left({\mathbf n}_1^{(j)},
{\mathbf n}_2^{(j)},{\mathbf n}_3^{(j)}=0\right)$$ is a non-mixing
sequence for $\alpha$ with ${\mathbf n}_s^{(j)}-{\mathbf
n}_t^{(j)}\to\infty$ as $j\to\infty$ for $s\neq t$. Then there are
elements $m_1,m_2,m_3$ of $R_2/\fp$, not all zero, with
\begin{equation}
\label{volochexample} m_1u^{{\mathbf n}_1^{(j)}}+ m_2u^{{\mathbf
n}_2^{(j)}}=-m_3
\end{equation}
for infinitely many $j$.
The field of fractions of $R_2/\fp$ may be identified
with ${\mathbb F}_2(t)$ by the map $u_1\mapsto t$,
$u_2\mapsto 1+t+t^2$, and in this field
(\ref{volochexample}) becomes
\begin{equation}\label{othervoloch}
ax+by=1
\end{equation}
with infinitely many solutions for $x,y$ in the
finitely generated multiplicative subgroup
$G=\langle\langle t,1+t+t^2\rangle\rangle$ of ${\mathbb F}_2(t)^{*}$.
By \cite{MR2000b:11029}, it follows that
(\ref{othervoloch}) is a {\sl$G$-trivial} equation:
there is an $n\ge 1$ for which $a^n,b^n\in G$.
Since $G$ is generated by irreducible polynomials,
this can only be true if $a,b\in G$. So there is an
infinite family of equations
\begin{equation}
\label{otherothervoloch} u^{{\mathbf m}_1^{(j)}}+u^{{\mathbf
m}_2^{(j)}} =1
\end{equation}
with ${\mathbf m}_1^{(j)}$, ${\mathbf m}_2^{(j)}$,
and
${\mathbf m}_1^{(j)}-{\mathbf m}_2^{(j)}\to\infty$ as $j\to\infty$.
By considering the shape of $\Newt(1+u_1+u_1^2+u_2)$,
this shows that the polynomial in (\ref{otherothervoloch})
has the same shape as $\Newt(1+u_1+u_1^2+u_2)$, so (without
loss of generality), ${\mathbf m}_1^{(j)}=(0,m(j))$
and ${\mathbf m}_2^{(j)}=(2m(j),0)$ for some $m(j)\to\infty$.
Thus the equation reduces to
\begin{equation}
\label{otherotherothervoloch}
(1+t+t^2)^{m(j)}=1+t^{2m(j)}.
\end{equation}
Write $m(j)=2^e\ell$, $\ell$ odd, for some $e\ge 0$.
Then the left-hand side of
(\ref{otherotherothervoloch}) is
\begin{eqnarray*}
(1+t+t^2)^{2^e\ell}&=&(1+t+O(t^2))^{2^e}\\
&=&1+t^{2^e}+O(t^2)^{2^e}\\
&=&1+t^{2^{e+1}\ell},
\end{eqnarray*}
which is impossible. It follows that $\mathscr M(\alpha)=3$.
\end{example}

\section{Further results}

\begin{remark}
(1) An extension of the arguments above using the product formula
for norms in function fields can be used to show that the set of
ratios of lengths of faces of $\Newt(A^{(j)})$ in a sequence that
witnesses the failure of $R$-fold mixing is bounded.

\noindent(2) If $S(f)$ is a rectangle, then the methods used
above show that the only non-mixing sequences of cardinality $4$
asymptotically arise from the non-mixing shape $S(f)$.
However, if $\vert S(f)\vert=4$ and $S(f)$ has a pair of
non-parallel sides, then this approach does not give
anything stronger than (1) above.
Since this approach does not give Conjecture \ref{conjr}
even for $r=4$ it has not been pursued further.

\noindent(3) Further progress on these problems, particularly
when $S(f)\cap\Newt(f)^{\circ}$ is non-empty, seems to require
Diophantine results on $S$-unit equations in finite characteristic.
\end{remark}

The method used to prove Proposition \ref{nastyprop} may also be
applied to prove the following theorem, the first part of which
relates to $\ZZ^d$-actions for any $d\ge2$.

\begin{theorem}{\rm(1)} If a sequence
$({\mathbf n}_1^{(j)},{\mathbf n}_2^{(j)},\dots,{\mathbf
n}_r^{(j)})$ in $\left(\ZZ^d\right)^r$ has the property that, for
every $s\neq t$,
$$
\frac{{\mathbf n}_s^{(j)}-{\mathbf n}_t^{(j)}}
{\Vert{\mathbf n}_s^{(j)}-{\mathbf n}_t^{(j)}\Vert}\longrightarrow
{\mathbf v}(s,t),
$$
for some vector ${\mathbf v}(s,t)$ whose entries are
linearly independent over $\mathbb Q$, then the sequence is
a mixing sequence for any mixing algebraic $\ZZ^d$-action.

\noindent{\rm(2)} Now fix $d=2$ and let $\alpha$ be determined by
an irreducible polynomial $f\in R_2^{(p)}$ as above. Then any
shape that does not contain all the faces of $\Newt(f)$ as
directions of differences is a mixing shape for $\alpha$.
\end{theorem}

\end{document}